\documentclass[a4paper,12pt]{amsart}

\usepackage{amsmath}
\usepackage{amssymb}
\usepackage{mathrsfs}
\usepackage{ifthen}
\usepackage{graphicx}
\usepackage[T1]{fontenc} 

\def\switchlinenumbers{\@ifstar
    {\let\makeLineNumberOdd\makeLineNumberRight
     \let\makeLineNumberEven\makeLineNumberLeft}%
    {\let\makeLineNumberOdd\makeLineNumberLeft
     \let\makeLineNumberEven\makeLineNumberRight}%
    }

\def\setmakelinenumbers#1{\@ifstar
  {\let\makeLineNumberRunning#1%
   \let\makeLineNumberOdd#1%
   \let\makeLineNumberEven#1}%
  {\ifx\c@linenumber\c@runninglinenumber
      \let\makeLineNumberRunning#1%
   \else
      \let\makeLineNumberOdd#1%
      \let\makeLineNumberEven#1%
   \fi}%
  }

\nonstopmode \numberwithin{equation}{section}
\newtheorem*{theorem*}{Theorem}

\newtheorem{thm}{Theorem}[section]
\newtheorem{cor}[equation]{Corollary}
\newtheorem{lem}[equation]{Lemma}
\newtheorem{prop}[equation]{Proposition}

\theoremstyle{definition}

\newtheorem{prob}[equation]{Problem}
\newtheorem{rem}{Remark}[section]

\newcounter{minutes}
\divide\time by 60
\newcounter{hours}
\multiply\time by 60
\addtocounter{minutes}{-\time}

\newcounter {own}
\def\theown {\thesection       .\arabic{own}}

\newenvironment{pf}[1][]{%
 \vskip 3mm
 \noindent
 \ifthenelse{\equal{#1}{}}%
  {{\slshape Proof. }}%
  {{\slshape #1.} }%
 }%
{\qed\bigskip}

\newcounter{alphabet}

\def\be{\begin{equation}}
\def\ee{\end{equation}}

\newcommand{\bee}{\begin{enumerate}}
\newcommand{\eee}{\end{enumerate}}

\newcommand{\blem}{\begin{lem}}
\newcommand{\elem}{\end{lem}}
\newcommand{\bthm}{\begin{thm}}
\newcommand{\ethm}{\end{thm}}
\newcommand{\bcor}{\begin{cor}}
\newcommand{\ecor}{\end{cor}}
\newcommand{\beg}{\begin{examp}}
\newcommand{\eeg}{\end{examp}}
\newcommand{\begs}{\begin{examples}}
\newcommand{\eegs}{\end{examples}}
\newcommand{\bdefe}{\begin{defin}}
\newcommand{\edefe}{\end{defin}}
\newcommand{\bprob}{\begin{prob}}
\newcommand{\eprob}{\end{prob}}
\newcommand{\bei}{\begin{itemize}}
\newcommand{\eei}{\end{itemize}}

\newcommand{\norm}[1]{\left\lVert#1\right\rVert}

\begin{document}

\title{Multidimensional Bohr radii for holomorphic functions with values in complex Banach spaces}

\author{Vasudevarao Allu}
\address{Vasudevarao Allu,
Department of Mathematics,
School of Basic Sciences,
Indian Institute of Technology Bhubaneswar,
Bhubaneswar-752050, Odisha, India.}
\email{avrao@iitbbs.ac.in}

\author{Himadri Halder}
\address{Himadri Halder,
Department of Mathematics, 
Indian Institute of Science, Bangalore-560012,
India}
\email{himadrihalder119@gmail.com}

\author{Subhadip Pal}
\address{Subhadip Pal,
	School of Mathematical Sciences,
	National Institute of Science Education and Research,
	Bhubaneswar, An OCC of Homi Bhabha National Institute, Jatni 752050,
	India.}
\email{subhadippal33@gmail.com}

\subjclass[{AMS} Subject Classification:]{Primary 32A05, 32A10, 46B45; Secondary 46G20}

\def\thefootnote{}
\footnotetext{ {\tiny File:~\jobname.tex,
printed: \number\year-\number\month-\number\day,
          \thehours.\ifnum\theminutes<10{0}\fi\theminutes }
} \makeatletter\def\thefootnote{\@arabic\c@footnote}\makeatother

\begin{abstract}

The main aim of this paper is to study multidimensional Bohr radii for holomorphic functions defined in complete Reinhardt domains in $\mathbb{C}^n$ with values in complex Banach spaces. More specifically, for holomorphic functions with values in arbitrary complex Banach spaces, we explore the asymptotic estimates of the classical Bohr radius and arithmetic Bohr radius in the unit ball of $\ell^n_q$ $(1\leq q\leq \infty)$ spaces. Further, we study a mixed version of Bohr radii for vector-valued holomorphic functions and as a consequence we obtain the exact value of mixed arithmetic Bohr radius.
  
\end{abstract}

\maketitle
\pagestyle{myheadings}
\markboth{Vasudevarao Allu, Himadri Halder, and Subhadip Pal}{Multidimensional Bohr radii for holomorphic functions with values in complex Banach spaces}

\section{Introduction}\label{section-1}
The Bohr radius problem plays a crucial role in connecting two major areas of research in pure mathematics, namely complex analysis and functional analysis. The Bohr radius problem's situation in multidimensional settings is much more complicated. Several research problems regarding the behavior of the multidimensional Bohr radius remain unresolved due to the lack of certain essential tools. The more interesting aspect of the Bohr radius problem is when we consider it in multidimensional settings for holomorphic functions with values in complex Banach spaces.

\vspace{2mm}

For the last two decades, a lot of attention has been paid to the generalizations of the multidimensional Bohr radius problem in various settings (see \cite{aizn-2000a,aizn-2005,boas-1997,defant-2003,defant-2007,popescu-2019}). One such generalization is studying Bohr's theorem for the power series of holomorphic functions defined in complete Reinhardt domains in $\mathbb{C}^n$. A complete Reinhardt domain $\Omega$ in $\mathbb{C}^n$ is a domain in $\mathbb{C}^n$ such that if $z=(z_{1},\ldots,z_{n}) \in \Omega$, then $(\eta_{1} e^{i\theta_{1}}z_{1},\ldots,\eta_{n} e^{i\theta_{n}}z_{n}) \in \Omega$ for all $\eta_{i} \in [0,1]$ and all $\theta_{i} \in \mathbb{R}$, $i=1, \ldots,n$. For a complex Banach space $X$, let $\mathcal{F}(\Omega , X)$ be the space of all holomorphic mappings $f$ in $\Omega$ into $X$ and $\mathcal{P}^k(\Omega,X)$ be the space of all $k$-homogeneous vector-valued polynomials, where $k \in \mathbb{N}$. Let $H_{\infty}(\Omega,X)$ be the set of all bounded holomorphic mappings on $\Omega$ into a complex Banach space $X$. In this article, we mostly focus on the aforesaid spaces with $\Omega=B_{\ell^n _q}$, where $q \in [1,\infty]$. For $q \in [1,\infty)$, 
$$B_{\ell^n _q}= \left\{z \in \mathbb{C}^n:\norm{z}_{q}=\left(\sum_{i=1}^{n}|z_{i}|^q\right)^{1/q}<1\right\}$$
 and $B_{\ell^n _\infty}= \{z \in \mathbb{C}^n: \norm{z}_{\infty}=\sup _{1\leq i \leq n} |z_{i}|<1\}$. For any given $p\in [1,\infty)$, $\lambda \geq 1$, and Reinhardt domain $\Omega$, the $\lambda_{p}$-Bohr radius $r^n_p(\Omega,X,\lambda)$ of $\Omega$ with respect to $\mathcal{F}(\Omega, X)$ is defined as the supremum of all $r \in [0,1]$ such that
 \begin{equation*} 
 	\sup_{z\in r\Omega}\sum_{\alpha\in \mathbb{N}^n_{0}} \norm{x_{\alpha}(f) z^{\alpha}}^p \leq \lambda \norm{f}^p _{\Omega}
  \end{equation*} 
 for all $f \in \mathcal{F}(\Omega, X)$ with $f(z)=\sum_{\alpha\in \mathbb{N}^n_{0}} x_{\alpha}(f) z^{\alpha}$ and $\norm{f}_{\Omega}=\sup \{\|f(z)\|: z \in \Omega\}$ and $\mathbb{N}_0=\mathbb{N}\cup\{0\}$. We write $r^n_p(\Omega,X,\lambda)=r_p(\Omega,X,\lambda)$ for $n=1$. \\

The remarkable result of Bohr \cite{Bohr-1914} states that $r_1(\mathbb{D},\mathbb{C},1)=1/3$. Much of the recent interest in Bohr's theorem stems from Dixon's work \cite{Dixon & BLMS & 1995}, when the author used it to disprove a long-standing open problem on the characterization of the Banach algebra satisfying the von Neumann inequality. Since then, a lot of research has been done to generalize Bohr's theorem in multidimensional and more abstract settings (see \cite{aizn-2000b}). A significant portion of this research focuses on studying Bohr's theorem in several variables for vector-valued functions, which led to the introduction of the quantity $r^n_p(\Omega,X,\lambda)$. However, determining the exact value of $r^n_p(\Omega,X,\lambda)$ is very challenging, and it seems there has not been much progress on this problem even for $X=\mathbb{C}$. For $\lambda \in [1,\sqrt{2}]$, Bombieri \cite{bombieri-1962} has shown that 
\begin{equation*} r_1(\mathbb{D},\mathbb{C},\lambda)=\frac{1}{3\lambda -2 \sqrt{2(\lambda^2 -1)}} \end{equation*} and for $\lambda$ very close to $\infty$, Bombieri and Bourgain \cite{bombieri-2004} have obtained the following estimate: \begin{equation*} r_1(\mathbb{D},\mathbb{C},\lambda) \sim \frac{\sqrt{\lambda^2 -1}}{\lambda}. \end{equation*}
In $1989$, Dineen and Timoney \cite{Dineen-Timoney-1989} extensively studied the constant $r^n_1(B_{\ell^n _\infty},\mathbb{C},1)$ and their result has been further investigated in \cite{boas-1997}. In $1997$, Boas and Khavinson \cite{boas-1997} obtained the following estimate for each $n \in \mathbb{N}$ with $n \geq 2$,
\begin{equation}\label{Pal-Vasu-P3-e-1.1}
	\frac{1}{3\sqrt{n}} < r^n_1(B_{\ell^n _\infty},\mathbb{C},1) < 2 \sqrt{\frac{\log n}{n}}.	
\end{equation}
The exact value of $r^n_1(B_{\ell^n _\infty},\mathbb{C},1)$ is still an open problem, and it has drawn the attention of many mathematicians to work further on the Bohr radius problem. Motivated by the above estimate of Boas and Khavinson, in $2000$, Aizenberg \cite{aizn-2000a} studied the constant $r^n_1(B_{\ell^n _q},\mathbb{C},1)$ with $q=1$ and proved that
\begin{equation}\label{Pal-Vasu-P3-e-1.2}   
	 \frac{1}{3e^{1/3}} < r^n_1(B_{\ell^n _1},\mathbb{C},1) \leq \frac{1}{3}.
\end{equation}

 In the same year, Boas \cite{boas-2000} generalized \eqref{Pal-Vasu-P3-e-1.1} and \eqref{Pal-Vasu-P3-e-1.2} to $r^n_1(B_{\ell^n _q},\mathbb{C},1)$ for $1 <q<\infty$. For $n>1$, Boas \cite{boas-2000} has proved that, if $1\leq q <2$, then 
\begin{equation} \label{itdn-e-1}
	\frac{1}{3\sqrt[3]{e}}\left(\frac{1}{n}\right)^{1-\frac{1}{q}} \leq r^n_1(B_{\ell^n _q},\mathbb{C},1) < 3 \left(\frac{\log \, n}{n}\right)^{1-\frac{1}{q}}
\end{equation}
and if $2 \leq q \leq \infty$, then
\begin{equation} \label{itdn-e-2}
	\frac{1}{3} \sqrt{\frac{1}{n}} \leq r^n_1(B_{\ell^n _q},\mathbb{C},1) < 2 \sqrt{\frac{\log\,n}{n}}.
\end{equation}
It is easy to see in \eqref{itdn-e-1} and \eqref{itdn-e-2} that the upper bounds contain a logarithmic factor that is not present in the lower bounds. Therefore, Boas observed that \eqref{itdn-e-1} and \eqref{itdn-e-2} do not give a sharp decay rate of Bohr radius with the dimension, and hence he conjectured that this logarithmic factor should not be really present. In $2006$, Defant {\it et al.} \cite{defant-2006} disproved this conjecture by obtaining a logarithmic lower bound that is almost correct asymptotic estimates for the Bohr radius $r^n_1(B_{\ell^n _q},\mathbb{C},1)$ with $1\leq q \leq \infty$. In particular, Defant {\it et al.} have proved that if $1\leq q \leq \infty$ then there is a constant $c>0$ such that
 \begin{equation} \label{itdn-e-3}    
 	\frac{1}{c} \left(\frac{\log\,n/\log\, \log\,n}{n}\right)^{1-\frac{1}{\min(q,2)}} \leq r^n_1(B_{\ell^n _q},\mathbb{C},1) \,\,\,\,\, \mbox{for all} \,\,\, n>1.
\end{equation}

In 2000, Djakov and Ramanujan \cite{Djakov & Ramanujan & J. Anal & 2000} first studied the Bohr radii $r_p(\mathbb{D}, \mathbb{C},1)$ and\\ $r^n_p(B_{\ell^n_{\infty}},\mathbb{C},1)$. For $p\in (1,2)$, the exact value of $r_p(\mathbb{D}, \mathbb{C},1)$ was determined by Kayumov and Ponnusamy \cite{Kay & Pon & AASFM & 2019}. In $2017$, Blasco \cite{Blasco-Collect-2017} extensively studied the constant $r^n_p(\Omega,X,\lambda)$ for $\Omega=\mathbb{D}, n=1, \lambda=1$ which is usually known as $p$-Bohr radius $r_p(\mathbb{D},X,1)$ for arbitrary Banach space $X$. It is worth noting that the precise value of $r_p(\mathbb{D},X,1)$ is very difficult to compute even for $X=\mathbb{C}$. Blasco \cite{Blasco-Collect-2017} has obtained a lower estimate of $r_p(\mathbb{D},\mathbb{C},1)$ and showed that $r_p(\mathbb{D},\mathbb{C},1) \geq p/(p+2)$ which is sharp only for $p=1$. 
\vspace{2mm}

Another aspect of the multidimensional Bohr radius is the arithmetic Bohr radius, which is more than a curious variant of the classical Bohr radius. In fact, we shall introduce new variants of the arithmetic Bohr radius, which we are going to define in the next two paragraphs.\\

Let $\Omega \subset \mathbb{C}^n$ be a complete Reinhardt domain. For each $1\leq p<\infty$ and $\lambda\geq 1$, the $\lambda_{p}$-\textit{arithmetic Bohr radius} of $\Omega$ with respect to $\mathcal{F}(\Omega, X)$ is defined as 
\begin{equation*}
	A_{p}(\mathcal{F}(\Omega, X), \lambda) := \sup \left\{\frac{1}{n}\sum_{i=1}^{n}r_i \,|\, r\in \mathbb{R}^{n}_{\geq 0},\, \forall\, f\in \mathcal{F}(\Omega, X) : \sum_{\alpha \in \mathbb{N}^{n}_{0}}\norm{x_{\alpha}(f)}^p r^{p\alpha} \leq \lambda \norm{f}^{p}_{\Omega}\right\},
\end{equation*} 
where $\mathbb{R}^n_{\geq 0}=\{r=(r_1,.\,.\,.\,, r_n) \in \mathbb{R}^n: r_i\geq 0, 1\leq i\leq n\}$ and $\norm{f}_{\Omega}=\sup \{\|f(z)\|: z \in \Omega\}$.  We write $A_p(\Omega,X,\lambda)$ for $A_p(H_{\infty}(\Omega,X),\lambda)$, $A_p(\Omega,X)$ for $A_p(\Omega,X,1)$ and $A^k_p(\Omega,X,\lambda)$ for $A_p(\mathcal{P}^k(\Omega,X),\lambda)$. For $X=\mathbb{C}$ and $p=1$, Defant {\it et al.} \cite{defant-2007} have first considered the constant $A_1(\mathcal{F}(\Omega,\mathbb{C}),\lambda)$ which is known as usual arithmetic Bohr radius. Recently, Kumar \cite{kumar-2023} has also studied the constant $A_1(\mathcal{F}(\Omega,\mathbb{C}),\lambda)$. The arithmetic Bohr radius $A_1(\mathcal{F}(\Omega,\mathbb{C}),\lambda)$ forms a useful tool to describe the domain of existence of the monomial expansion of complex-valued holomorphic functions in a complete Reinhardt domain (cf. \cite{prengel-2005}). This led us to introduce the notion of $\lambda_{p}$-arithmetic Bohr radius $A_{p}(\mathcal{F}(\Omega, X), \lambda)$ for arbitrary complex Banach spaces, which plays a crucial role in estimating sequences in the set of monomial convergence of vector-valued holomorphic functions in a Reinhardt domain $\Omega$ in $\mathbb{C}^n$. One of the main purposes of this paper is to investigate the constant $A_{p}(\mathcal{F}(\Omega, X), \lambda)$ by establishing a relation between $\lambda_{p}$- Bohr radius and $\lambda_{p}$-arithmetic Bohr radius. Our results mainly generalize the results obtained in \cite{defant-2007}.

\vspace{2mm}
We now concentrate on another multidimensional generalization of the Bohr radius $R_{p,q}(X)$ which has been introduced for the first time by Blasco \cite{Blasco-OTAA-2010} in order to study the Bohr inequality for Banach space valued holomorphic functions in the unit disk $\mathbb{D}$. In this paper, we introduce a new notion of $R_{p,q}(X)$ in more general multidimensional settings. Let $\Omega \subset \mathbb{C}^n$ be a complete Reinhardt domain. For each $1\leq p, q<\infty$ and $\lambda\geq 1$, we define $\lambda_{p,q}$-\textit{Bohr radius}, denoted by $R^n_{p,q}(\mathcal{F}(\Omega, X), \lambda)$, of $\Omega$ with respect to $\mathcal{F}(\Omega, X)$ to be the supremum of all $r \in [0,1]$ such that
\begin{equation}\label{Pal-Vasu-P3-e-1.5}
\sup_{z\in r\Omega}	\left(\norm{x_0}^p+
	\left(\sum_{m=1}^{\infty}\sum_{|\alpha|=m}\norm{x_{\alpha}(f)z^{\alpha}} \right)^q \right) \leq \lambda \norm{f}_{\Omega}
\end{equation}
for all $f \in \mathcal{F}(\Omega, X)$ with $f(z)=\sum_{k=0}^{\infty} \sum_{|\alpha|=k} x_{\alpha}(f) z^{\alpha}$. We use the convention $R^n_{p,q}(\Omega,X,\lambda)$ instead of $R^n_{p,q}(\mathcal{F}(\Omega, X),\lambda)$ and $R^n_{p,q}(\Omega,X,\lambda)=R_{p,q}(\Omega,X,\lambda)$ for $n=1$. 
Blasco \cite{Blasco-OTAA-2010} has extensively studied the Bohr radius $R_{p,q}(\mathbb{D},X,1)$ for complex Banach spaces $X$. In the same paper, it has been established \cite[Theorem 1.2]{Blasco-OTAA-2010} that the classical Bohr radius $R(\mathbb{D}, \ell^n_p,1)=0$ for $1\leq p<\infty$. To make this modified Bohr radius more significant, we introduce the constant $\lambda$ in \eqref{Pal-Vasu-P3-e-1.5}. Determining the exact value of $\lambda_{p,q}$-Bohr radius is difficult in-general for different complex Banach spaces. Recently, Das \cite{das-2023} has obtained the exact value of $R_{p,q}(\mathbb{D},\mathbb{C},1)$ for $p,q\in[1,\infty)$ and the asymptotic behavior of $R^n_{p,q}(B_{\ell^n_{\infty}},X,1)$ for infinite dimensional Banach spaces. In this situation, our aim is to study $\lambda_{p,q}$-Bohr radius $R_{p,q}(\mathbb{D},X,\lambda)$ in terms of another notion of arithmetic Bohr radius, namely, $\lambda_{p,q}$-\textit{arithmetic Bohr radius} $A_{p,q}(\mathcal{F}(\Omega, X), \lambda)$ of $\Omega$ with respect to $\mathcal{F}(\Omega,X)$, which is defined by 
\begin{align*}
	&A_{p,q}(\mathcal{F}(\Omega, X), \lambda) :=\\& \sup \left\{\frac{1}{n}\sum_{i=1}^{n}r_i \,|\, r\in \mathbb{R}^{n}_{\geq 0},\, \mbox{for all}\, f\in \mathcal{F}(\Omega, X):\norm{x_0}^p+
	\left(\sum_{m=1}^{\infty}\sum_{|\alpha|=m}\norm{x_{\alpha}(f)} r^{\alpha}\right)^q \leq \lambda \norm{f}_{\Omega}\right\},
\end{align*}
where $\mathbb{R}^n_{\geq 0}=\{r=(r_1,.\,.\,., r_n) \in \mathbb{R}^n: r_i\geq 0, 1\leq i\leq n\}$ and $\norm{f}_{\Omega}=\sup \{\|f(z)\|: z \in \Omega\}$.\\[2mm]

The organization of this paper is as follows: In Theorem \ref{Pal-Vasu-P3-thm-1}, we provide the estimate of $\lambda_p$-arithmetic Bohr radius of bounded holomorphic functions in terms of $\lambda_p$-arithmetic Bohr radius of $k$-homogeneous polynomials in complete Reinhardt domains. Theorem \ref{Pal-Vasu-P3-thm-2.2} gives the asymptotic estimates of the $\lambda_p$-arithmetic Bohr radius for the $n$-dimensional $\ell^n_q$-balls. Next, Theorem \ref{Pal-Vasu-P3-prop-2.4} relates the classical $\lambda_{p}$-Bohr radius and the $\lambda_p$-arithmetic Bohr radius in the unit ball of $\ell^n_q$ spaces. We obtain an important estimate of the arithmetic Bohr radius of homogeneous polynomials in Proposition \ref{Pal-Vasu-P3-prop-2.4}. We establish the connection between $\lambda_{p,q}$-arithmetic Bohr radii of different vector-valued holomorphic function spaces and the bounded holomorphic function spaces in Theorem \ref{Pal-Vasu-P3-thm-2.6}. Finally, we investigate the exact value of $\lambda_{p,q}$-arithmetic Bohr radius of unit ball in $\ell^n_1$ space through the proof of Theorem \ref{Pal-Vasu-P3-thm-2.7}.

\section{Main Results}

Our first result provides an estimate of $\lambda_p$-arithmetic Bohr radius of bounded holomorphic functions in terms of $\lambda_p$-arithmetic Bohr radius of $k$-homogeneous polynomials in complete Reinhardt domains. 
\begin{thm}\label{Pal-Vasu-P3-thm-1}
	Let $\Omega$ be a complete Reinhardt domain in $\mathbb{C}^n$, $\lambda\geq 1$, $1\leq p <\infty$, and $C=\max \{1/3, ((\lambda-1)/(2\lambda-1))^{1/p}\}$. Then we have 
	\begin{equation*}
		C \,A_p\left(\bigcup_{k=1}^{\infty}\mathcal{P}^{k}(\Omega),\mathbb{C}, \lambda\right) \leq A_p\left(H_{\infty}(\Omega),\mathbb{C}, \lambda\right) \leq A_p\left(\bigcup_{k=1}^{\infty}\mathcal{P}^k(\Omega), \mathbb{C},\lambda\right).
	\end{equation*}  
\end{thm}

\noindent 
It is worth mentioning that the result of Defant {\it et al.} \cite[Proposition 2.2]{defant-2007} follows from Theorem \ref{Pal-Vasu-P3-thm-1} for $p=1$. In the following result, Theorem \ref{Pal-Vasu-P3-thm-2.2}, we are going to give the asymptotic estimates of the $\lambda_p$-arithmetic Bohr radius for the $n$-dimensional $\ell^n_p$-balls. The main ingredients of our proof of Theorem \ref{Pal-Vasu-P3-thm-2.2} are the following facts, which we record as Lemma \ref{Pal-Vasu-P3-lem-2.2}, Lemma \ref{P3-boas-lem-2.3}, and Proposition \ref{Pal-Vasu-P3-prop-2.4}. To make the statement of these results brief, let us recall the following convention from \cite{defant-2007}.
For bounded Reinhardt domains $\Omega_{1}, \Omega_{2}\subset \mathbb{C}^n$, let
\begin{equation*}
	S(\Omega_1, \Omega_{2}):= \inf \left\{s>0 : \Omega_{1} \subset s\Omega_{2}\right\}.
\end{equation*}
 By a Banach sequence space $Y$, we mean a complex Banach space $Y\subset \mathbb{C}^{\mathbb{N}}$ such that $\ell_{1}\subset Y \subset \ell_{\infty}$. If $\Omega$ is a bounded Reinhardt domain in $\mathbb{C}^n$ and $Y$ and $Z$ are Banach sequence spaces, we write
\begin{equation}\label{Pal-Vasu-P3-e-2.1}
	S(\Omega, B_{Y_n})= \sup_{z\in \Omega}\norm{z}_{Y} \quad \mbox{and }\quad  S(B_{Y_n}, B_{Z_n})= \norm{\mbox{id}: Y_n \rightarrow Z_n},
\end{equation}
where $Y_n$(resp. $Z_n$) is the space spanned by first $n$ canonical basis vectors $e_n$ in $Y$(resp. $Z$). Throughout the paper, $X$ will stand for a complex Banach space. The next result compares the $\lambda_p$-arithmetic Bohr radii for two bounded Reinhardt domains in $\mathbb{C}^n$.
\begin{lem}\label{Pal-Vasu-P3-lem-2.1}
	Let $\Omega_{1}, \Omega_{2}$ be bounded Reinhardt domains in $\mathbb{C}^n$ and $1\leq p<\infty$. Then for $\lambda \geq 1$, we have 
	\begin{equation*}
		A_p(\Omega_{1},X, \lambda) \leq S(\Omega_{1}, \Omega_{2})\, A_p(\Omega_{2}, X, \lambda).
	\end{equation*} 
\end{lem}

\begin{rem}
	\begin{enumerate}
		\item[(a)] For a bounded Reinhardt domain $\Omega$ in $\mathbb{C}^n$, it is easy to observe that $S(\Omega, t\Omega)=1/t$ and $S(t\Omega, \Omega)=t$  for all $t>0.$ Therefore, an immediate consequence from Lemma \ref{Pal-Vasu-P3-lem-2.1} we have that for any $1\leq p<\infty$,
		\begin{equation*}
			A_p(t\Omega_{1}, X, \lambda)=tA_p(\Omega_{1}, X, \lambda),\quad \mbox{for all}\quad t>0.
 		\end{equation*}
 	\item[(b)]$A_p(\cdot,X, \lambda)$ is increasing, that is, $A_p(\Omega_{1}, X,  \lambda)\leq A_p(\Omega_{2},X, \lambda)$ whenever $\Omega_{1} \subset \Omega_{2}.$
 	
 	\item[(c)] It is obvious that $A_p(\Omega,X,\lambda)\leq A^k_p(\Omega,X,\lambda)$.
 	\item[(d)] We note that Lemma \ref{Pal-Vasu-P3-lem-2.1} generalizes the result of Defant {\it et. al.} \cite[Lemma 4.2]{defant-2007}.
	\end{enumerate}
	
\end{rem}

Using the above notation, the following lemma gives a connection between the $\lambda_p$-arithmetic Bohr radius $A_p(\Omega, X, \lambda)$ for bounded Reinhardt domain $\Omega\subset \mathbb{C}^n$ and the $p$-Bohr radius $r^n_p(\Omega,X,\lambda)$ of a complex Banach space $X$, which helps us to obtain the lower bound in Theorem \ref{Pal-Vasu-P3-thm-2.2}.
\begin{lem}\label{Pal-Vasu-P3-lem-2.2}
	Let $\Omega \subset\mathbb{C}^n$ be a bounded Reinhardt domain and $1\leq p<\infty$. Then for $\lambda \geq 1$ we have,
	\begin{equation*}
		A_p(\Omega,X, \lambda) \geq \frac{S(\Omega, B_{\ell^n_1})}{n}r^n_p(\Omega,X,\lambda).
	\end{equation*}
\end{lem}

We need the following estimate for $k$-homogeneous polynomials $(k\geq 2)$ in $\mathbb{C}^n$ which has been proved by Boas \cite{boas-2000}.

\begin{lem}\cite{boas-2000}\label{P3-boas-lem-2.3}
	Let $1\leq q \leq \infty$ and $n$ and $k$ are positive integers greater than $1$. Then there exists a homogeneous polynomial of degree $k$ in $\mathbb{C}^n$ and signs $\epsilon_{\alpha}\in \{-1, 1\}$ such that 
	\begin{equation*}
		\bigg\|\sum_{|\alpha|=k}\epsilon_{\alpha}\frac{k!}{\alpha!}z^{\alpha}\bigg\|_{B_{\ell^n_q}}\leq \sqrt{32k\log (6k)}\,\times\, \left\{\begin{array}{ll}
			n^{\frac{1}{2}}(k!)^{1-\frac{1}{q}}  & \mbox{ if  \,$q\leq 2$},\\[1mm]
			n^{\frac{1}{2}+\left(\frac{1}{2}-\frac{1}{q}\right)k}(k!)^{\frac{1}{2}}  & \mbox{ if \,$q\geq 2$}.
		\end{array}\right.
	\end{equation*}
\end{lem}

Now we obtain an upper bound for $\lambda_p$-arithmetic Bohr radius of $k$-homogeneous polynomials in $\mathbb{C}^n$ which later gives an upper bound for $\lambda_{p}$-arithmetic Bohr radius of $B_{\ell^n_q}\,\,(1\leq q\leq \infty)$ with respect to $H_{\infty}(B_{\ell^n_q}, X).$
\begin{prop}\label{Pal-Vasu-P3-prop-2.4}
	Let $1\leq p< \infty$ and $\lambda \geq 1$. Then for all $1\leq q\leq \infty$ and $k, n \geq 2$ we have 
	\begin{equation*}
		A^k_p(B_{\ell^n_q},X, \lambda) \leq d \sqrt[pk]{\lambda} \sqrt[2k]{k\log k (k!)^{2(1-(1/\min\{2, q\}))}n}\, \left(\frac{1}{n}\right)^{\frac{1}{p}+\frac{1}{\max\{2, q\}}-\frac{1}{2}},
	\end{equation*}
where $d$ is a uniform constant.
\end{prop}

 Following the techniques from \cite{defant-2007, defant-2004} and all above preparation, we now wish to state one of our main results with regard to the asymptotic estimates of the $\lambda_p$-arithmetic Bohr radius for the $n$-dimensional $\ell^n_q$-balls.
\begin{thm}\label{Pal-Vasu-P3-thm-2.2}
	Let $1\leq p< \infty$ and $\lambda \geq 1$. Then there exists a constant $d>0$ such that for all $1\leq q\leq \infty$ and $n>1$, we have 
	\begin{equation}\label{Pal-Vasu-P3-e-2.8}
		\frac{r^n_p(B_{\ell^n_q},X,\lambda)}{n^{1/q}}\leq A_p(B_{\ell^n_q},X, \lambda)\leq d\, \frac{\left(\log n\right)^{1-\left(1/\min\{q, 2\}\right)}}{n^{\frac{1}{p}+\frac{1}{\max\{2, q\}}-\frac{1}{2} }}.
	\end{equation}
\end{thm}
\begin{rem}
	In  Theorem \ref{Pal-Vasu-P3-thm-2.2}, for $p=1$ we obtain the upper estimate of \cite[Theorem 4.1]{defant-2007}, but we feel that the lower estimate of \eqref{Pal-Vasu-P3-e-2.8} can be improved for arbitrary complex Banach spaces.
\end{rem}

Next, we present a relationship between arithmetic Bohr radius and geometric property of complex Banach spaces. Recently, Das \cite{das-2023} has studied the $p$-Bohr radius $r_p(B_{\ell^n_{\infty}}, X,1)$ of a complex Banach space $X$. In the same paper \cite[Theorem 1.5]{das-2023}, the author has established the asymptotic estimates of $r_p(B_{\ell^n_{\infty}}, X,1)$ for any $p$-uniformly $PL$-convex $(p\geq 2)$ complex Banach spaces $X$. A complex Banach space $X$ is said to be $p$-uniformly $PL$-convex $(2\leq p<\infty)$ if there exists a constant $\delta>0$ such that
\begin{equation}\label{Pal-Vasu-P3-e-2.11}
	\norm{x}^p+\delta \norm{y}^p\leq \frac{1}{2\pi}\int_{0}^{2\pi}\norm{x+e^{i\theta}y}^p\,d\theta
\end{equation}
for all $x, y \in X$. We denote $I_p(X)$ be the supremum of all such $\delta>0$ satisfying \eqref{Pal-Vasu-P3-e-2.11}.\\[1mm]

In view of \cite[Theorem 1.5]{das-2023}, Lemma \ref{Pal-Vasu-P3-lem-2.2} and the fact that $S(B_{\ell^n_{\infty}}, B_{\ell^n_1})=1/n$, we have the following consequence.
\begin{cor}
	Let $B_{\ell^n_{\infty}}$ be the unit polydisk in $\mathbb{C}^n$ and $2\leq p<\infty$. Then for all $\lambda \geq 1$ and $n>1$, we have
	\begin{equation*}
		A_p(B_{\ell^n_{\infty}},X, \lambda)\geq \frac{1}{n^2} \left(\frac{I_p(X)}{2^p+I_p(X)}\right)^{2/p},
	\end{equation*}
where $X$ is any $p$-uniformly PL-convexity complex Banach space $X$ with dim$(X)\geq 2$.
\end{cor}

Now we are in position to give an interesting connection between the classical $\lambda_{p}$-Bohr radius and  $\lambda_p$-arithmetic Bohr radius in unit ball of $\ell^n_q$ spaces.

\begin{thm}\label{Pal-Vasu-P3-thm-2.4}
	Let $1\leq \lambda <\infty$ and $1\leq p<\infty$. Then for every $n$ and $1\leq q<\infty$ we have 
	\begin{equation*}
		\frac{r_p(\mathbb{D},X,\lambda)}{n}\leq A_p(B_{\ell^n_q},X,\lambda)\leq \left(\frac{r_p(\mathbb{D},X,\lambda)}{n^{1/p}}\right)^{1/q}.
	\end{equation*}
\end{thm}

We want to recall that we have provided the notion of $\lambda_{p,q}$-arithmetic Bohr radius in Section \ref{section-1}. In the following main result, we establish that the $\lambda_{p,q}$-arithmetic Bohr radii of different vector-valued holomorphic function spaces $\mathcal{O}$ corresponds to $\lambda_{p,q}$-arithmetic Bohr radius of bounded holomorphic function space if  $\mathcal{O}$ contains all vector-valued polynomials.
\begin{thm}\label{Pal-Vasu-P3-thm-2.6}
	Let $\Omega$ be a complete Reinhardt domain in $\mathbb{C}^n$ and $\lambda\geq 1$. If $\mathcal{O}(\Omega, X)$ is a space of all holomorphic mappings on $\Omega$ into a complex Banach space $X$ containing all vector-valued polynomials then for every $1\leq p,q, <\infty$, we have
	\begin{equation*}
		A_{p,q}(\mathcal{O}(\Omega, X), \lambda)=A_{p,q}(H_{\infty}(\Omega, X), \lambda).
	\end{equation*}
\end{thm}
An analogous version of Theorem \ref{Pal-Vasu-P3-thm-2.4} can also be shown in the setting of $\lambda_{p,q}$-arithmetic Bohr radius which later provides us an exact value of arithmetic Bohr radius.
\begin{thm}\label{Pal-Vasu-P3-thm-2.7}
	Let $1\leq q < \infty$ and $\lambda \geq 1$. Then for every $1\leq p_{1}, q_{1}< \infty$ and $n$ we have 
	\begin{equation*}
		\frac{R_{p_{1}, q_{1}}(\mathbb{D}, X,\lambda)}{n}\leq A_{p_1,q_1}(B_{\ell^n_{q}}, X,\lambda) \leq \left(\frac{R_{p_{1}, q_{1}}(\mathbb{D}, X,\lambda)}{n}\right)^{1/q}. 
	\end{equation*}
\end{thm}

In 2010, O. Blasco \cite[Proposition 1.4]{Blasco-OTAA-2010}, established the exact value of Bohr radius $R_{p,1}(\mathbb{D},\mathbb{C})$ as $p/(p+2)$ for $1\leq p\leq 2$ and $\lambda=1$. In view of of Blasco's result, we record the following easy consequence of Theorem \ref{Pal-Vasu-P3-thm-2.7} as the exact value of arithmetic Bohr radius of the unit ball in $\ell^n_{1}$.

\begin{cor}
	For every $1\leq p\leq 2$ and $n\in \mathbb{N}$, we have $A_{p,1}(B_{\ell^n_{1}},\mathbb{C},1)=p/(n(p+2)).$
\end{cor}
\section{Proof of Main Results}
\begin{pf} [{\bf Proof of Theorem \ref{Pal-Vasu-P3-thm-1}}]
	The right-hand inequality is clear from the fact that 
	\begin{equation*}
		\bigcup_{k=1}^{\infty} \mathcal{P}^{k}(\Omega) \subset H_{\infty} (\Omega).
	\end{equation*}
	Choose $r\in \mathbb{R}^n_{\geq 0}$ such that for all $k\geq 1$ and for all $k$-homogeneous polynomial $h \in \mathcal{P}^k(\mathbb{C}^n)$, we have
	\begin{equation}\label{Pal-Vasu-P3-e-3.1}
		\sum_{|\alpha|=k} |c_{\alpha}(h)r^{\alpha}|^p \leq \lambda \norm{h}^p_{\Omega}.
	\end{equation}
	We wish to show that 
	\begin{equation}\label{Pal-Vasu-P3-e-3.2}
		\max \left\{\frac{1}{3}, \left(\frac{\lambda-1}{2\lambda-1}\right)^{1/p}\right\}\frac{1}{n}\sum_{i=1}^{n}r_i \leq A_p(H_{\infty}(\Omega)_,\, \lambda).
	\end{equation}
	Fix $f \in H_{\infty}(\Omega)$. Without loss of generality, we assume that $\norm{f}_{\infty}=1$ and consider its monomial series expansion as $f(z)=\sum_{\alpha\in \mathbb{N}^n_{0}}c_{\alpha}(f)z^{\alpha}$.\\[2mm]
	Given $z\in \Omega$, we define the function $g: \mathbb{D}\rightarrow \mathbb{C}$ by
	\begin{equation*}
		g(\zeta)=f(\zeta z)=\sum_{\alpha\in \mathbb{N}^n_{0}}c_{\alpha}(f)\zeta^{|\alpha|}z^{\alpha}=\sum_{k=0}^{\infty}\left(\sum_{|\alpha|=k}c_{\alpha}(f)z^{\alpha}\right)\zeta^{k}.
	\end{equation*}
	Let $\theta \in \mathbb{R}$ be such that $e^{i\theta}c_{0}(f)=|c_0(f)|$ and consider the function $F(\zeta)=1-e^{i\theta}g(\zeta)$, $\eta \in \mathbb{D}$. It is worth to notice that $\zeta z\in \Omega$ for every $\zeta \in \mathbb{D}$, thus we have $\norm{g}_{\infty}\leq 1$ on $\mathbb{D}$ and 
	\begin{equation*}
		\mbox{Re} \, F(\zeta)=\mbox{Re} \left(1-e^{i\theta}g(\zeta)\right)\geq 0.
	\end{equation*}
	Suppose $F(\zeta)=\sum_{k=0}^{\infty}a_k\zeta^k$, then 
	\begin{equation*}
		\sum_{k=0}^{\infty}a_k\zeta^k=1-\sum_{k=0}^{\infty}\left(e^{i\theta}\sum_{|\alpha|=k}c_{\alpha}(f)z^{\alpha}\right)\zeta^k=1-e^{i\theta}c_0(f)+ \sum_{k=1}^{\infty}\left(e^{i\theta}\sum_{|\alpha|=k}c_{\alpha}(f)z^{\alpha}\right)\zeta^k.
	\end{equation*}
	It follows that $a_0=1-|c_0(f)|$ and $a_k=e^{i\theta}\sum_{|\alpha|=k}c_{\alpha}(f)z^{\alpha}$, for $k\in \mathbb{N}.$\\[2mm]
	In view of Carath\'eodory's theorem \cite{aizn-2005}, we obtain that 
	\begin{equation}\label{Pal-Vasu-P3-e-3.3}
		\left|\sum_{|\alpha|=k}c_{\alpha}(f)z^{\alpha}\right|\leq 2(1-|c_0(f)|).
	\end{equation}
	From \eqref{Pal-Vasu-P3-e-3.1} and \eqref{Pal-Vasu-P3-e-3.3}, we have
	\begin{align*}
		\sum_{\alpha\in \mathbb{N}^n_{0}}\left|c_{\alpha}(f)\left(\frac{r}{3}\right)^{\alpha}\right|^p & =
		|c_0(f)|^p+\sum_{k=1}^{\infty}\sum_{|\alpha|=k}\left|c_{\alpha}(f)\left(\frac{r}{3}\right)^{\alpha}\right|^p\\ & = |c_0(f)|^p +\sum_{k=1}^{\infty} \frac{1}{3^{pk}} \sum_{|\alpha|=k} |c_{\alpha}(f)|^p r^{p\alpha}\\ &\leq |c_0(f)|^p + \sum_{k=1}^{\infty} \frac{1}{3^{pk}}\, \lambda \sup_{z\in \Omega} \left|\sum_{|\alpha|=k}c_{\alpha}(f)z^{\alpha}\right|^p\\ & \leq 
		|c_0(f)|^p+ \lambda\, \frac{2^p}{3^p-1}\left(1-|c_0(f)|\right)^p \leq \lambda=\lambda \norm{f}_{\Omega}.
	\end{align*}
	This gives one of the estimates in \eqref{Pal-Vasu-P3-e-3.2}. To obtain the other estimate, we set $s=({(2\lambda-1)}/{(\lambda-1)})^{1/p}.$ By going similar lines of argument as above and using the well-known Cauchy's inequalities, we obtain
	
	\begin{align*}
		\sum_{\alpha\in \mathbb{N}^n_{0}}\left|c_{\alpha}(f)\left(\frac{r}{s}\right)^{\alpha}\right|^p &= |c_0(f)|^p +\sum_{k=1}^{\infty} \frac{1}{s^{pk}} \sum_{|\alpha|=k} |c_{\alpha}(f)|^p r^{p\alpha}\\ & \leq |c_0(f)|^p + \lambda \sum_{k=1}^{\infty} \frac{1}{s^{pk}} \norm{P_k(f)}^p_{\Omega}\\ &\leq 
		|c_0(f)|^p+\lambda \norm{f}^p_{\Omega} \sum_{k=1}^{\infty}\frac{1}{s^{pk}} \leq \left(1+\frac{\lambda}{s^p-1}\right)\norm{f}^p_{\Omega}=\lambda \norm{f}^p_{\Omega}.
	\end{align*}
	Therefore, $(1/sn)\sum_{i=1}^{n}r_i\leq A_p(H_{\infty}(\Omega), \mathbb{C},\lambda).$ This completes the proof.
\end{pf}

\begin{pf} [{\bf Proof of Lemma \ref{Pal-Vasu-P3-lem-2.1}}]
	Let $r=(r_1, \cdots, r_n)\in \mathbb{R}^n_{\geq 0}$ be such that 
	\begin{equation*}
		\sum_{\alpha\in \mathbb{N}^n_{0}} \norm{x_{\alpha}(f)r^{\alpha}}^p\leq \lambda
		\norm{f}^p_{\Omega_1}
	\end{equation*}
	for all $f(z)=\sum_{\alpha\in \mathbb{N}^n_{0}}x_{\alpha}(f)z^{\alpha}\in H_{\infty}(\Omega_{1}, X)$. Set $t:= S(\Omega_{1}, \Omega_{2})$. Then for all $\epsilon >0$, we have 
	
	\begin{equation*}
		\sum_{\alpha\in \mathbb{N}^n_{0}}\norm{x_{\alpha}(f)}^p\left(\frac{r}{t+\epsilon}\right)^{p\alpha
		}\leq \lambda \bigg\|\sum_{\alpha\in \mathbb{N}^n_{0}}x_{\alpha}(f)\left(\frac{z}{t+\epsilon}\right)^{\alpha}\bigg\|^p_{\Omega_{1}} \leq \lambda \norm{f}^p_{\Omega_{2}}
	\end{equation*}
	for all $f\in H_{\infty}(\Omega_{2}, X).$ By letting $\epsilon \rightarrow 0$, we obtain 
	\begin{equation*}
		\frac{1}{nt}\sum_{i=1}^{n}r_i \leq A_p(\Omega_{2}, X, \lambda).
	\end{equation*}
	Therefore, it follows that $	A_p(\Omega_{1},X, \lambda) \leq S(\Omega_{1}, \Omega_{2})\, A_p(\Omega_{2}, X, \lambda).$
\end{pf}
\begin{pf} [{\bf Proof of Lemma \ref{Pal-Vasu-P3-lem-2.2}}]
In view of \eqref{Pal-Vasu-P3-e-2.1}, we have 
	\begin{equation*}
		S(\Omega, B_{\ell^n_{1}})=\sup_{z\in \Omega}\norm{z}_{\ell^n_{1}}.
	\end{equation*}
	Therefore, for given $0<\epsilon<r^n_p(\Omega,X,\lambda)$, we can find an element $\tilde{z}\in \Omega$ such that 
	\begin{equation*}
		\norm{\tilde{z}}_{\ell^n_1}\geq S(\Omega, B_{\ell^n_{1}})-\epsilon.
	\end{equation*}
	Let $s:=r^n_p(\Omega,X,\lambda)-\epsilon$, $w:=s\tilde{z}$, and $r:=s|\tilde{z}|=|w|.$ Since $w\in s\Omega$ and $s<r^n_p(\Omega,X,\lambda)$, for $f=\sum_{\alpha\in \mathbb{N}^n_{0}}x_{\alpha}z^{\alpha}\in H_{\infty}(\Omega, X)$, we have
	\begin{equation*}
		\sum_{\alpha\in \mathbb{N}^n_{0}}\norm{x_{\alpha}}^p r^{p\alpha}=\sum_{\alpha\in \mathbb{N}^n_{0}}\norm{x_{\alpha}w^{\alpha}}^p\leq \bigg\|\sum_{\alpha\in \mathbb{N}^n_{0}}x_{\alpha}z^{\alpha}\bigg\|^p_{s\Omega}\leq \lambda\norm{f}^p_{\Omega}.
	\end{equation*}
	Therefore, we obtain 
	\begin{equation*}
		A_p(\Omega, X,\lambda) \geq \frac{1}{n}\sum_{i=1}^{n}r_i = \frac{r^n_p(\Omega,X,\lambda)-\epsilon}{n}\norm{\tilde{z}}_{\ell^n_{1}}\geq \frac{r^n_p(\Omega,X,\lambda)-\epsilon}{n}\left(S(\Omega, B_{\ell^n_{1}})-\epsilon\right)
	\end{equation*}
	holds for all $\epsilon>0.$ This completes the proof.
\end{pf}

\begin{pf} [{\bf Proof of Proposition \ref{Pal-Vasu-P3-prop-2.4}}]
	Let $r=(r_1,.\,.\,., r_n)\in \mathbb{R}^n_{\geq 0}$ be such that 
	\begin{equation*}
		\sum_{|\alpha|=k}\norm{x_{\alpha}(f)}^pr^{p\alpha} \leq \lambda \norm{f}^p_{B_{\ell^n_q}}
	\end{equation*}
	for all $f\in \mathcal{P}^k(\mathbb{C}^n, X)$. From H\"older's inequality, we have $\norm{r}^p_{1}\leq n^{p-1}\norm{r}^p_{p}$. In view of Lemma \ref{P3-boas-lem-2.3}, there exist signs $\epsilon_{\alpha}\in\{-1,1\}$ such that
	\begin{equation*}
		\left(\frac{1}{n}\sum_{i=1}^{n}r_i\right)^{pk}\leq
	\frac{1}{n^k}\left(\sum_{i=1}^{n}r^p_{i} \right)^k	= \frac{1}{n^k}\sum_{|\alpha|=k}\frac{k!}{\alpha!}r^{p\alpha}\leq \frac{\lambda}{n^k}\bigg\|\sum_{|\alpha|=k}\epsilon_{\alpha}\frac{k!}{\alpha!}z^{\alpha}\bigg\|^p_{B_{\ell^n_q}}\leq \frac{\lambda}{n^k}R^p_k,
	\end{equation*}
	where $R_k$ is the expression in the right-hand side of Lemma \ref{P3-boas-lem-2.3}. Taking the $pk$-th root of both sides, we obtain
	\begin{equation*}
		A^k(B_{\ell^n_q},X, \lambda) \leq \frac{\sqrt[pk]{\lambda}\sqrt[k]{R_k}}{n^{1/p}}.
	\end{equation*}
	This completes the proof.
\end{pf}
\vspace{-4mm}
\begin{pf} [{\bf Proof of Theorem \ref{Pal-Vasu-P3-thm-2.2}}]
	In view of Lemma \ref{Pal-Vasu-P3-lem-2.2} and using the fact that 
	\begin{equation*}
		S(B_{\ell^n_q}, B_{\ell^n_1})=n^{1-(1/q)} \quad\, \mbox{for all}\,\,\,1\leq q\leq \infty,
	\end{equation*}
	we obtain the lower bound of \eqref{Pal-Vasu-P3-e-2.8}. For the upper bound of $A_p(B_{\ell^n_q}, X, \lambda)$, by Proposition \ref{Pal-Vasu-P3-prop-2.4}, it is enough to show that there exists a constant $d>0$ such that for all $n>1$ we are able to find a $k\,(\geq 2)$ with
	\begin{equation}\label{Pal-Vasu-P3-e-2.9}
		\sqrt[pk]{\lambda} \sqrt[2k]{k \log k (k!)^{2\left(1-\left(1/\min \{2, q\}\right)\right)}n} \leq d\, \left(\log n\right)^{1-\left(1/\min \{2, q\}\right)}.
	\end{equation}
	It is worth to note that the quantity $\sqrt[2k]{k\log k}$ is bounded in $k$ and also $\sqrt[k]{k!}\leq k/\sqrt{2}.$ So we can estimate the left-hand side of \eqref{Pal-Vasu-P3-e-2.9} upto a constant independent of $k$ and $n$ by
	\begin{equation}\label{Pal-Vasu-P3-e-2.10}
		\sqrt[pk]{\lambda}\,k^{1-\left(1/\min\{2, q\}\right)}\sqrt[2k]{n}.
	\end{equation}
	Let us choose $k=2$ for $n=2,.\,.\,.,8$ and $k=[\log n]$ (the greatest integer less than or equal to $\log n$) for $n \geq 9$. Then, the above expression \eqref{Pal-Vasu-P3-e-2.10} is bounded above by 
	\begin{equation*}
		\lambda^{\frac{2}{p\log n}}(\log n)^{1-\frac{1}{\min \{2, q\}}}e, 
	\end{equation*}
	which gives the desired bound. This completes the proof.
\end{pf}
\vspace{-2mm}
\begin{pf} [{\bf Proof of Theorem \ref{Pal-Vasu-P3-thm-2.4}}]
	Let $r=r_p(\mathbb{D},X,\lambda)$ and $f(z)=\sum_{\alpha\in \mathbb{N}^n_{0}}x_{\alpha}(f)z^{\alpha} \in H_{\infty}(B_{\ell^n_q},X)$.  Consider the function $g(z)=f(\xi z)$, where $\xi=(1, 0,.\,.\,.,0)\in B_{\ell^n_q}$ and $z\in \mathbb{D}$. Clearly, $g\in H_{\infty}(\mathbb{D},X)$ with $\norm{g}_{\mathbb{D}}=\norm{f}_{B_{\ell^n_q}}$.
	Therefore, we have 
	\begin{equation*}
		\sum_{\alpha\in \mathbb{N}^n_{0}}\norm{x_{\alpha}(f)}^p (r,0,\,.\,.\,.,0)^{p\alpha}=\sum_{k=0}^{\infty}\norm{x_{k\xi}(f)}^pr^{pk}=\sum_{k=0}^{\infty}\norm{x_k(g)}^pr^{pk}\leq \lambda \norm{g}^p_{\mathbb{D}}=\lambda\norm{f}^p_{B_{\ell^n_q}}.
	\end{equation*}
	This shows that $(r/n)\leq A_p(B_{\ell^n_q},X,\lambda)$, and hence we obtain $(r_p(\mathbb{D},X,\lambda)/n) \leq A_p(B_{\ell^n_q},X,\lambda).$  Conversely, we prove that 
	\begin{equation*}
		A_p(B_{\ell^n_q},X,\lambda)\leq \left(\frac{r_p(\mathbb{D},X,\lambda)}{n^{1/p}}\right)^{1/q}.
	\end{equation*}
	Suppose $r\in \mathbb{R}^n_{\geq 0}$ such that for all $h \in H_{\infty}(B_{\ell^n_q},X)$,
	\begin{equation*}
		\sum_{\alpha\in \mathbb{N}^n_{0}}\norm{x_{\alpha}(h)}^pr^{p\alpha}\leq \lambda\norm{h}^p_{B_{\ell^n_q}}.
	\end{equation*}
	We want to prove that 
	\begin{equation*}
		\frac{\sum_{j=1}^{n}r_j}{n}\leq \left(\frac{r_p(\mathbb{D},X,\lambda)}{n^{1/p}}\right)^{1/q}.
	\end{equation*}
	Let $f:\mathbb{D}\rightarrow X$ be a bounded holomorphic map and we consider the function $s: B_{\ell^n_q}\rightarrow \mathbb{D}$ defined by 
	\begin{equation*}
		s(z)=z^q_1+\cdots+z^q_n,\quad z\in B_{\ell^n_q}.
	\end{equation*}
	Now if we set $h=f\circ s$, then we have $\norm{h}_{B_{\ell^n_q}}=\norm{f}_{\mathbb{D}}.$
	Moreover, for each $z\in B_{\ell^n_{q}}$,
	\begin{equation*}
		h(z)=\sum_{k=0}^{\infty}x_k(f)s(z)^k=\sum_{k=0}^{\infty}x_k(f)\sum_{|\alpha|=k}\frac{k!}{\alpha!}z^{q\alpha}=\sum_{\alpha\in \mathbb{N}^n_{0}}x_{\alpha}(h)z^{q\alpha},
	\end{equation*}
	where 
	\begin{equation*}
		x_{\alpha}(h)=\frac{k!}{\alpha!}x_k(f)
	\end{equation*}
	whenever $|\alpha|=k$. Then for all $z\in \mathbb{C}^n$, we have
	\begin{align*}
		\sum_{\alpha\in \mathbb{N}^n_{0}}\norm{x_{\alpha}(h)}^p|z|^{pq\alpha}&= \sum_{k=0}^{\infty}\norm{x_k(f)}^p\sum_{|\alpha|=k}\left(\frac{k!}{\alpha!}\right)^p|z|^{pq\alpha}\\[2mm]& \geq \sum_{k=0}^{\infty} \norm{x_k(f)}^p\sum_{|\alpha|=k}\frac{k!}{\alpha!}|z|^{pq\alpha} =\sum_{k=0}^{\infty}\norm{x_k(f)}^p\norm{z}^{pqk}_{pq}.
	\end{align*}
	An observation shows that
	\begin{equation*}
		\sum_{k=0}^{\infty}\norm{x_k(f)}^p\norm{r}^{pqk}_{pq}\leq \sum_{\alpha\in \mathbb{N}^n_{0}}\norm{x_{\alpha}(h)r^{q\alpha}}^p \leq \lambda \norm{h}^p_{B_{\ell^n_{q}}}=\lambda\norm{f}^p_{\mathbb{D}},
	\end{equation*}
	which implies that $\norm{r}^q_{pq}\leq r_p(\mathbb{D},X,\lambda).$ Again we have $\norm{r}^{pq}_{1}\leq n^{(pq-1)}\norm{r}^{pq}_{pq}.$ Hence we obtain $n^{(1-pq)}\norm{r}^{pq}_{1}\leq r^p_p(\mathbb{D},X,\lambda)$. This completes the proof.
\end{pf}

\begin{pf} [{\bf Proof of Theorem \ref{Pal-Vasu-P3-thm-2.6}}]
	We start with an observation that 
	\begin{equation*}
		\mathcal{P}(\Omega, X)\subseteq \mathcal{O}(\Omega, X)\subseteq H_{\infty}(\Omega, X).
	\end{equation*}
	Then, we have the following inclusion
	\begin{equation*}
		A_{p,q}(H_{\infty}(\Omega, X), \lambda) \leq A_{p,q}(\mathcal{O}(\Omega, X), \lambda) \leq A_{p,q}(\mathcal{P}(\Omega, X), \lambda).
	\end{equation*}
	It suffices to show that for any $\lambda\geq 1$, $A_{p,q}(\mathcal{P}(\Omega, X), \lambda)\leq 	A_{p,q}(H_{\infty}(\Omega, X), \lambda)$.\\[2mm]
	Assume $r\in \mathbb{R}^n_{\geq 0}$ such that for all $P\in \mathcal{P}(\mathbb{C}^n, X)$,
	\begin{equation}\label{Pal-Vasu-P3-e-2.14}
		\norm{x_0(P)}^p+\left(\sum_{k=1}^{\infty}\norm{x_{\alpha}(P)}r^{\alpha}\right)^q\leq \lambda \norm{P}_{\Omega}.
	\end{equation}
	Fix $f\in H_{\infty}(\Omega, X)$ and $0<t<1$. Our aim is to show that 
	\begin{equation}\label{Pal-Vasu-P3-e-2.15}
		\norm{x_0(f)}^p+\left(\sum_{k=1}^{\infty}\norm{x_{\alpha}(f)}(tr)^{\alpha}\right)^q\leq \lambda \norm{f}_{\Omega}.
	\end{equation}
	For $k\geq 0$, consider the $k$-homogeneous polynomials $\chi_{k}=\sum_{|\alpha|=k}x_{\alpha}(f)z^{\alpha}$, where $x_{\alpha}(f)$ is the $\alpha$-th coefficient in the monomial series expansion $\sum_{\alpha\in \mathbb{N}^n_{0}}x_{\alpha}(f)z^{\alpha}$ of $f$. Note that $\Omega$ is balanced, so by Cauchy inequalities we have 
	\begin{equation}\label{Pal-Vasu-P3-e-2.16}
		\norm{\chi_{k}}_{t\Omega}=t^k\norm{\chi_{k}}_{\Omega}\leq t^k \norm{f}_{\Omega}.
	\end{equation}
	Therefore, the series $\sum_{k=0}^{\infty}\chi_{k}$ converges absolutely and uniformly to $f$ on $t\Omega$. Thus, given $\epsilon>0$, there exits $n_0 \in \mathbb{N}$ such that for all $n \geq n_0$,
	\begin{equation*}
		\bigg\|\sum_{k=0}^{n}\chi_{k}\bigg\|_{t\Omega}\leq \norm{f}_{t\Omega} + \epsilon \leq \norm{f}_{\Omega}+\epsilon.
	\end{equation*}
	We now consider the polynomials $\Gamma_{k}(z)=\chi_{k}(tz)$ for $z\in \mathbb{C}^n$. Then,
	\begin{equation*}
		x_{\alpha}(\Gamma_{k})=
		\left\{\begin{array}{ll}
			t^kx_{\alpha}(f) & \mbox{ if  \,$|\alpha|=k$},\\[1mm]
			0  & \mbox{ if \,$|\alpha|\neq k$}.
		\end{array}\right.
	\end{equation*}
	From \eqref{Pal-Vasu-P3-e-2.14} and \eqref{Pal-Vasu-P3-e-2.16}, we obtain that
	\begin{align*}
		\norm{x_0(f)}^p+\left(\sum_{k=1}^{\infty}\norm{x_{\alpha}(f)}(tr)^{\alpha}\right)^q &= \norm{x_0(\Gamma_{0})}^p + \left(\sum_{\alpha\in \mathbb{N}^n}\bigg\|x_{\alpha}\left(\sum_{k=1}^{n}\Gamma_{k}\right)r^{\alpha}\bigg\|\right)^q\\[2mm] & \leq \lambda \bigg\|\sum_{k=0}^{n}\Gamma_{k}\bigg\|_{\Omega}=\lambda \bigg\|\sum_{k=0}^{n}\chi_{k}\bigg\|_{t\Omega}\leq \lambda \left(\norm{f}_{\Omega}+\epsilon\right),
	\end{align*} 
	which holds for all $\epsilon>0$. Consequently, this proves \eqref{Pal-Vasu-P3-e-2.15}.
	Hence for every $r\in \mathbb{R}^n_{\geq 0}$ satisfying \eqref{Pal-Vasu-P3-e-2.14} and for all $f\in H_{\infty}(\Omega, X)$ and $0<t<1$, we have shown that 
	\begin{equation*}
		t\sum_{i=1}^{n}r_i \leq A_{p,q}(H_{\infty}(\Omega,X),\lambda),
	\end{equation*}
	which proves our claim. This completes the proof.
\end{pf}
\vspace{-2mm}

\begin{pf} [{\bf Proof of Theorem \ref{Pal-Vasu-P3-thm-2.7}}]
	First we show the left-hand inequality $(R_{p_1, q_1}(\mathbb{D}, X, \lambda))/n\leq A_{p_1,q_1}(B_{\ell^n_q},X, \lambda).$ Assume $r=R_{p_1,q_1}(\mathbb{D},X,\lambda)$ and $f\in H_{\infty}(B_{\ell^n_q},X)$. We define $g(z)=f(ze_1)=f(z,0,.\,.\,.,0)$ for $z\in \mathbb{D}.$ Then $g\in H_{\infty}(\mathbb{D},X)$ and $\norm{g}_{\mathbb{D}}=\norm{f}_{B_{\ell^n_q}}$. Further,
	\begin{align*}
		\norm{x_0(f)}^{p_1} + \left(\sum_{k=1}^{\infty}\norm{x_{\alpha}(f)}(r,0,.\,.\,.,0)^{\alpha}\right)^{q_1} &= 	\norm{x_0(f)}^{p_1} + \left(\sum_{k=1}^{\infty}\norm{x_{ke_1}(f)}r^k\right)^{q_1}\\&= \norm{x_0(g)}^{p_1} + \left(\sum_{k=1}^{\infty}\norm{x_k(g)}r^k\right)^{q_1}\\[2mm] & \leq \lambda \norm{g}_{\mathbb{D}}=\lambda \norm{f}_{B_{\ell^n_q}}.
	\end{align*}
	Therefore, we have $r/n=A_{p_1,q_1}(B_{\ell^n_q},X, \lambda)$, which gives the desired inequality. \\[2mm]
	On the other hand, we want to prove that $	A_{p_1,q_1}(B_{\ell^n_q},X, \lambda) \leq \left({R_{p_1,q_1}(\mathbb{D},X,\lambda)}/{n}\right)^{1/q}.$
	Let $r\in \mathbb{R}^n_{\geq 0}$ be such that 
	\begin{equation*}
		\norm{x_0(u)}^{p_1} + \left(\sum_{k=1}^{\infty}\sum_{|\alpha|=k}\norm{x_{\alpha}(u)}r^{\alpha}\right)^{q_1}\leq \lambda \norm{u}_{B_{\ell^n_q}}
	\end{equation*}
	for all $u\in H_{\infty}(B_{\ell^n_q},X)$. Now our aim is to show that 
	\begin{equation*}
		\frac{1}{n}\left(\sum_{j=1}^{n}r_j\right) \leq \left(\frac{R_{p_1,q_1}(\mathbb{D},X,\lambda)}{n}\right)^{1/q}.
	\end{equation*}
	Fix $f\in H_{\infty}(\mathbb{D},X)$, and define the function
	\begin{equation*}
		v(z)=z^q_1+\cdots+z^q_n, \quad z=(z_1,.\,.\,.,z_n)\in B_{\ell^n_q}. 
	\end{equation*}
	Then for $u=f\circ v$, we have $\norm{u}_{B_{\ell^n_q}}=\norm{f}_{\mathbb{D}}$. Also for each $z\in B_{\ell^n_q}$, we obtain
	\begin{equation*}
		u(z)=\sum_{k=0}^{\infty}x_{k}(f)v(z)^k=\sum_{k=0}^{\infty}x_k(f)\sum_{|\alpha|=k}\frac{k!}{\alpha!}z^{q\alpha}=\sum_{\alpha\in \mathbb{N}^n_{0}}x_{\alpha}(u)z^{q\alpha},
	\end{equation*}
	where $x_{\alpha}(u)=x_k(f)(k!/\alpha!)$ whenever $|\alpha|=k$. Then for all $z\in \mathbb{C}^n$, we have
	\begin{align*}
		\norm{x_0(u)}^{p_1} + \left(\sum_{k=1}^{\infty}\norm{x_{\alpha}(u)}|z|^{q\alpha}\right)^{q_1}&= 
		\norm{x_0(f)}^{p_1} + \left(\sum_{k=1}^{\infty}\norm{x_k(f)}\sum_{|\alpha|=k}\frac{k!}{\alpha!}|z|^{q\alpha}\right)^{q_1}\\[2mm] & =\norm{x_0(f)}^{p_1} + \left(\sum_{k=1}^{\infty}\norm{x_k(f)}\norm{z}^{qk}_{q}\right)^{q_1},
	\end{align*}
	and hence
	\begin{align*}
		\norm{x_0(f)}^{p_1} + \left(\sum_{k=1}^{\infty}\norm{x_k(f)}\norm{r}^{qk}_{q}\right)^{q_1} & =
		\norm{x_0(u)}^{p_1} + \left(\sum_{k=1}^{\infty}\norm{x_{\alpha}(u)}r^{q\alpha}\right)^{q_1}\\[2mm]&\leq \lambda \norm{u}_{B_{\ell^n_q}}=\lambda \norm{f}_{\mathbb{D}}.
	\end{align*}
	Therefore, it follows that $\norm{r}^q_{q}\leq R_{p_1,q_1}(\mathbb{D},X,\lambda)$. Since we have $\norm{r}^q_{1}\leq n^{q-1}\norm{r}^q_{q}$, finally we obtain $n^{1-q}\norm{r}^q_{1}\leq R_{p_1,q_1}(\mathbb{D},X,\lambda)$. This completes the proof. 
\end{pf}

 \vspace{3mm}
 
 \noindent\textbf{Compliance of Ethical Standards:}\\
 
 \noindent\textbf{Conflict of interest.} The authors declare that there is no conflict  of interest regarding the publication of this paper.
 \vspace{1.5mm}
 
 \noindent\textbf{Data availability statement.}  Data sharing is not applicable to this article as no datasets were generated or analyzed during the current study.
 \vspace{1.5mm}
 
 \noindent\textbf{Authors contributions.} Both the authors have made equal contributions in reading, writing, and preparing the manuscript.\\
 
\noindent\textbf{Acknowledgment:} 
 The authors would like to express their sincerest gratitude to the referee for careful reading of the manuscript and many valuable suggestions, which greatly helped to improve the clarity of the exposition in this manuscript. 
The research of the first named author was supported by SERB-CRG, Govt. of India.  The research of the second named author was supported by Institute Post-Doctoral Fellowship of IIT Bombay, and the research of the third named author was supported by DST-INSPIRE Fellowship (IF 190721),  New Delhi, India.

\end{document}